\documentclass[10.9pt,twoside]{amsart}
\usepackage{amsmath, amsthm, amscd, amsfonts, amssymb, graphicx, color}
\usepackage[bookmarksnumbered, colorlinks, plainpages]{hyperref}

\theoremstyle{definition}

\theoremstyle{remark}

\numberwithin{equation}{section}

\begin{document}
\setcounter{page}{1}
\begin{center}
{\bf BOUNDED TOPOLOGICAL GROUPS }
\end{center}

\title[]{}
\author[]{KAZEM HAGHNEJAD AZAR  }

\address{}

\dedicatory{}



\begin{abstract} In this note  for a topological group $G$, we introduce a bounded subset of $G$ and we find some relationships of this definition with other topological properties of $G$.

\end{abstract} \maketitle

\section{\bf  preliminaries and
Introduction }

\noindent  Suppose that $G$ is a topological group and $E\subseteq G$. In this paper, we want to know when $E$ is bounded or unbounded subset of $G$ and if $G$ is metrizablity, we show that $E\subseteq G$ is bounded with respect to topology if and only if it is bounded with respect to metric. Let $E\subseteq G$ be bounded and closed. Then $E$ is compact subset of $G$. Conversely if $E$ is a component of $e$ and compact, then $E$ is bounded. We investigated some topological property for bounded subset of $G$.\\
Now we introduce some notations and definitions that we used
throughout  this paper.\\
For topological group $G$, $e$ is identity element of $G$ and for $E\subseteq G$, $E^-$ is closure of $E$ and for every $n\in\mathbb{N}$,
$$E^n=\{x_1x_2x_3...x_n:~~x_i\in E, 1\leq i\leq n\}.$$
A topological space $X$ is $O-dimensional$ if the family of all sets that are both open and closed is open basis for the topology.\\\\

\begin{center}
\section{ \bf Bounded Topological Groups  }
\end{center}

\noindent{\it{\bf Definition 2-1.}} Let $G$ be topological group and $E\subseteq G$. We say that $E$ is bounded subset of $G$, if for every neighborhood $V$ of $e$, there is natural number $n$ such that $E\subseteq V^n$.\\\\
 \noindent It is clear that if $E$ is bounded subset of $G$ and $H$ is subgroup of $G$, then $E/H$ is bounded subset of $G/H$.\\\\

\noindent{\it{\bf Theorem 2-2.}} Let $G$ be topological group and metrizable with respect to a left invariant metric $d$. Then $G$ is bounded with respect to topology if and only if $G$ is bounded with respect to metric $d$.\\

\begin{proof} Let $G$ be a bounded topological group and $\varepsilon>0$. Take $d([0,\varepsilon))=U\times V$ where $U$ and $V$ are neighborhoods of $e$. Suppose that $W$ is symmetric  neighborhood of $e$ such that $W\subseteq U\cap V$. Then there is natural number $n$ such that $W^n=G$. Since $d(W\times W)<\varepsilon$, we show that $d(W^2\times W^2)<2\varepsilon$, and so $d(W^n\times W^n)<n\varepsilon$.
Assume that $x,y,x^\prime, y^\prime\in W$. Then we have
$$d(xy, x^\prime y^\prime)\leq d(xy, e)+d(e,x^\prime y^\prime)=d(y, x^{-1})+d({x^\prime}^{-1}, y^\prime)<2\varepsilon.$$
Then $d(G\times G)=d(W^n\times W^n)<n\varepsilon$.\\
Conversely, suppose that $G$ is bounded with respect to metric $d$. Then there is $M>0$ such that $d(G\times G)<M$. Let $U$ be a neighborhood of $e$. Choose $\varepsilon >0$ such that $d^{-1}([0,\varepsilon))\subseteq U\times U$. Thak natural number $n$ such that $n\varepsilon >M$. Then we have
$$G\times G=d^{-1}([0,M))=d^{-1}([0,n\varepsilon))\subseteq V^n\times V^n.$$
It follows that $G=V^n$, and so that $G$ is bounded.

\end{proof}

\noindent{\it{\bf Theorem 2-3.}} Let $G$ be topological group and let $H$ be a normal subgroup of $G$. If $H$ and $G/H$ are bounded, then $G$ is bounded.

\begin{proof} Let $U$ be a neighborhood of $e$. Put $V=U\cap H$. Then there are natural numbers $m$ and $n$ such that

$(U/H)^n=G/H$ and $V^m=H$. We show that $U^{n+m}=G$.\\
Let $x\in G$. Then if $x\in H$, we have
$$x\in V^m\subset U^m\subset U^{n+m}.$$
Now let $x\notin H$. Then $xH\in (U/H)^n$. Assume that $x_1, x_2, ...,x_n\in U$ such that
$$xH=x_1x_2...x_nH.$$
Consequently there is $h\in H$ such that $xh\in U^n$, and so $x\in U^nH\subset U^nV^m\subset U^nU^m=U^{n+m}$. We conclude that $U^{n+m}=G$, and so $G$ is bounded.
\end{proof}

\noindent{\it{\bf Theorem 2-4.}} If $G$ is a locally compact O-dimensional topological group, then $G$ is unbounded.
\begin{proof} Let $U$ be a neighborhood of $e$ such that $U^-$ is compact and $U^-\neq G$. Since $G$ is a O-dimensional topological group, $U$ contains an open and closed neighborhood as $V$. Then $V$ is a compact neighborhood of $e$. By apply [1, Theorem 4.10] to obtain a neighborhood $W$ of $e$ such that $WV\subset V$. Take $W_0=W\cap V$. Then ${W_0}^2\subset WV\subset V\subset U^-$. By finite induction, we have
$${W_0}^n\subset W_0{W_0}^{n-1}\subset WV\subset V\subset U^-,$$
for every natural number $n$. It follows that ${W_0}^n\subsetneqq G$ for every natural number $n$, and so $G$ is unbounded.
\end{proof}

\noindent{\it{\bf Theorem 2-5.}} Suppose that  $G$ is a locally compact, Hausdorff, and totally disconnected  topological group. Then $G$ is unbounded.
\begin{proof}
By using [1, Theorem 3.5] and Theorem 2-4, proof is hold.

\end{proof}

\noindent{\it{\bf Theorem 2-6.}} Let $G$ be topological group. Then we have the following assertions.
\begin{enumerate}
\item  If $E\subseteq G$ is bounded, then $E^-$ is bounded subset of $G$.
\item  If $G$ is bounded, then $G$ is connected and moreover $G$ has no proper open subgroups.
\end{enumerate}
\begin{proof} 1) Let $U$ be a neighborhood of $e$ and suppose that $V$ is a neighborhood of $e$ such that $V^-\subset U$. Since $E$ is bounded subset of $G$, there is natural number $n$ such that $E\subset V^n$. Then $E^-\subset {(V^n)}^-\subset {(V^-)}^n\subset U^n$. It follows that $E^-$ is a bounded subset of $G$.\\
2) Since $G$ is bounded, there is a natural number $n$ such that $G=V^n$ where $V$ is neighborhood of $e$. By using [1, Corollary 7.9], proof is hold.

\end{proof}

\noindent{\it{\bf Corollary 2-7.}} Assume that $G$ is a locally compact  topological group. Then every bounded and closed subset of $G$ is compact, moreover if $E\subseteq G$  is bounded, then $E^-$ is compact.\\\\
Every bounded topological group $G$, in general, is not compact, for example $\mathbb{R}/\mathbb{Z}$ is bounded, but is not compact.\\\\

\noindent{\it{\bf Theorem 2-8.}}  Let $G$ be topological group and suppose that $E\subseteq G$  is the component of $e$. If $E$ is compact, then $E$ is bounded.
\begin{proof} Since $E$ is the component of $e$, by using [1, Theorem 7.4], for every neighborhood $U$ of $e$, we have $E\subseteq \bigcup_{k=1}^\infty U^k$. Since $E$ is compact there is natural number $n$ such that $E\subseteq U^n$. Then $E$ is bounded subset of $G$.
\end{proof}

In general, every compact subset $E$ of a topological group $G$ is not bounded and in above Theorem, it is necessary that $E$ must be a component of $e$. For example $Z_n=\{ \bar{0}, \bar{1}, \bar{2},...,\bar{n}\}$ for every $n\geq 1$, with discrete topology is not bounded, but it is compact.\\\\
\noindent{\it{\bf Corollary 2-9.}} If $G$ is a locally compact  topological group, then the component of $e$ is bounded.\\\\
\noindent{\it{\bf Theorem 2-10.}}  Let $G$ and $G^\prime$ be topological group and suppose that  $\pi :G\rightarrow G^\prime$ is group isomorphism. If $\pi$ is continuous and $E\subseteq G$ is a bounded subset of $G$, then $\pi(E)$ is bounded subset of $G^\prime$.
\begin{proof} Let $V^\prime$ be a neighborhood of $e^\prime\in G^\prime$. Then $\pi^{-1}(V^\prime)$ is a neighborhood of $e$. Since $E$ is a bounded subset of $G$, there is a natural number $n$ such that $E\subseteq (\pi^{-1})^n(V^\prime)\subseteq \pi^{-1}({V^\prime}^n)$ implies that $\pi(E)\subseteq {V^\prime}^n$. Thus $\pi(E)$ is a bounded subset of $G^\prime$.
\end{proof}

\noindent{\it{\bf Definition 2-11.}} Let $G$ and $G^\prime$ be topological group. We say that the mapping   $\pi :G\rightarrow G^\prime$ is compact, if for every bounded subset $E\subseteq G$, $\pi(E)$ is relatively compact.\\\\

\noindent{\it{\bf Theorem 2-12.}} Let $G$ and $G^\prime$ be topological group and suppose that  $\pi :G\rightarrow G^\prime$ is continuous and  group isomorphism. Then if $G^\prime$ is locally compact, then $\pi$ is compact.
\begin{proof} Let $E\subseteq G$ be bounded. By using Theorem 2.10, $\pi(E)$ is bounded subset of $G^\prime$ and by using Theorem 2.6, $\pi(E)^-$ is compact, and so that $\pi$ is compact.

\end{proof}

\bibliographystyle{amsplain}

\end{document}